\documentclass[11pt]{amsart}
\usepackage[paperwidth=17cm, paperheight=22.5cm, bottom=2.5cm, right=2.5cm]{geometry}
\usepackage{amssymb,amsmath,amsthm} 
\usepackage[english]{babel}
\usepackage[utf8x]{inputenc}
\usepackage[T1]{fontenc}
\usepackage{enumerate}
\usepackage{mathtools}
\usepackage{graphicx}
\usepackage{subfigure} 
\usepackage{float}
\graphicspath{{Img/}} 
\usepackage[export]{adjustbox}
\usepackage{tikz}
\usetikzlibrary{matrix}
\usepackage{multirow}

\newtheorem{teo}{Theorem}[section]
\newtheorem{propo}[teo]{Proposition}
\newtheorem{cor}[teo]{Corollary}
\newtheorem{lem}[teo]{Lemma}
\newtheorem{que}{Question}
\newtheorem*{que*}{Question}

\theoremstyle{definition}
\newtheorem{dfn}[teo]{Definition}

\theoremstyle{remark}

\newtheorem{fact}[teo]{Fact}

\newcommand{\rest}{\upharpoonright} 
\newcommand{\forces}{\Vdash}		

\newcommand{\baire}{\omega^{\omega}}

\def\cA{{\mathcal{A}}} \def\cB{{\mathcal{B}}} \def\cC{{\mathcal{C}}} \def\cD{{\mathcal{D}}} \def\cE{{\mathcal{E}}} \def\cF{{\mathcal{F}}}  \def\cH{{\mathcal{H}}} \def\cI{{\mathcal{I}}}  \def\cK{{\mathcal{K}}}     \def\cP{{\mathcal{P}}}   \def\cS{{\mathcal{S}}}

\title{Madness and (weak) normality}
\author{C\'esar Corral}
\address{Centro de Ciencias  Matem\'aticas, Universidad Nacional Aut\'onoma de M\'exico,
\' Campus Morelia, 58089, Morelia, Michoac\'an, M\'exico.}
\email{cicorral@matmor.unam.mx}

\date{}

\thanks{%
The author gratefully acknowledges support from CONACyT scholarship 742627.}

\begin{document}

\begin{abstract}
We consider weakenings of normality in $\Psi$-spaces and prove that the existence of a MAD family whose $\Psi$-space is almost-normal is independent of \textsf{ZFC}. We also construct a partly-normal not quasi-normal AD family, answering  questions of Garc\'ia-Balan and Szeptycki. We finish by showing that the concepts of almost-normal and strongly $\aleph_0$-separated AD families are different, even under \textsf{CH}, answering a question of Oliveira-Rodrigues and Santos-Ronchim.
\end{abstract}

\maketitle

\section{Introduction and notation}

Two subsets $A,B\subseteq\omega$ are \emph{almost disjoint} if $|A\cap B|<\omega$. A family $\cA\subseteq\cP(\omega)$ is \emph{almost disjoint} (AD for short) if its elements are pairwise almost disjoint. We say that $\cA$ is \emph{maximal almost disjoint} (MAD) if it is AD and maximal with respect to this property (equivalently, for every infinite $X\subseteq\omega$, there exists $A\in\cA$ such that $|A\cap X|=\omega$).

Each AD family $\cA$ has a topological space naturally associated with it. The \emph{$\Psi$-space} or the \emph{Mr\'owka-Isbell space} associated to $\cA$ is denoted by $\Psi(\cA)$ and the underlying set is $\omega\cup\cA$, where the points in $\omega$ are isolated and for every $A\in\cA$, sets of the form $\{A\}\cup (A\setminus F)$ are basic neighborhoods for $A$, with $F$ ranging over all finite subsets of $\omega$.

It follows easily from the definition that $\Psi(\cA)$ is a separable, locally compact, zero dimensional, scattered, first countable Moore space (see \cite{mrowka1955completely} and \cite{tall1977set}). Despite their simplicity, almost disjoint families and their $\Psi$-spaces are central tools in set theoretic topology, since many important problems have an equivalent reformulation in the realm of $\Psi$-spaces. In particular, every separable hereditarily locally compact space is homeomorphic to a $\Psi$-space \cite{kannan1979hereditarily}.

Mro\'owka-Isbell spaces provide a wide and numerous source of examples and counterexamples in many areas of topology. Many examples of the use of AD families and their $\Psi$-spaces can be found in \cite{hruvsak2014almost}. Normality is no exception. 
A MAD family is never normal, AD families of size $\mathfrak{c}$ are not normal by Jones' lemma, since $\cA$ is a discrete subspace of size continuum of a separable space. One of the first examples of an AD family with special combinatorial properties, was a Luzin family \cite{luzin1947subsets}. An AD family $\cA$ is a \emph{Luzin family}, if it can be enumerated as $\cA=\{A_\alpha:\alpha<\omega_1\}$ in such a way that $\{\beta<\alpha:A_\alpha\cap A_\beta\subseteq n\}$ is finite for every $\alpha<\omega_1$ and every $n\in\omega$. The key property of Luzin families is that if $\cB,\cC\subseteq\cA$ are two uncountable subfamilies, they can not be separated, in consequence, Luzin families are not normal. This suggest that normality is not easily fulfilled for an AD family. 
One of the first applications of AD families to problems related to normality, was the equivalence of the existence of a normal, separable, non-metrizable Moore space and the existence of an uncountable AD family which is not normal. The later was proved to be independent of \textsf{ZFC} \cite{tall1977set}.

In \cite{szeptycki2020weak}, weak normality properties on $\Psi$-spaces were considered.
Recall that a space $X$ is \emph{normal} if every two disjoint closed sets $C,D\subseteq X$ can be separated by two disjoint open sets $U,V\subseteq X$ (that is $C\subseteq U$, $D\subseteq V$ and $U\cap V=\emptyset$).  A subset $C\subseteq X$ of a topological space is \emph{regular closed} if $C=\overline{int(C)}$. Thus, the definition of normality becomes weaker if we require one, or both of the closed sets to be regular closed or a finite intersection of regular closed sets (which is called $\pi$-closed). Ranging over these possibilities, several weakenings of normality arise, and so do some implications between them (see, \cite{alshammari2019partial}, \cite{alshammari2020quasi}, \cite{singal1970almost} and \cite{szeptycki2020weak}). We summarize these implications in the next diagram without defining all the concepts involved simply to organize them and have a visual support. We will define each term that we will focus on, when necessary.
\begin{itemize}
    \item[($\ast$)] $\textnormal{normal}\implies\textnormal{almost-normal}\implies \textnormal{quasi-normal}\implies$
    \item[] $$\textnormal{partly-normal}\implies \textnormal{mildly-normal.}$$
\end{itemize}

Counterexamples of some of these implications were given in \cite{szeptycki2020weak}: A mildly-normal which is not partly-normal and a quasi-normal which is not almost-normal AD families were constructed, whilst counterexamples of the remaining two implications were left open. In particular, the existence of an almost-normal MAD family, was left open (Questions 4.1, 4.2 and 4.3 in\cite{szeptycki2020weak}). 
In Section \ref{section2}, we provide an example of an almost-normal MAD family under \textsf{CH}. In Section \ref{section3} we show that under \textsf{PFA}, no MAD family can be almost-normal, proving that the existence of an almost-normal MAD family is independent of the axioms of \textsf{ZFC}. In Section \ref{section4}, we build a partly-normal AD family which is not quasi-normal, hence, completing all the counterexamples in $(\ast)$, at least, consistently. Finally, in Section \ref{section5}, we will construct a strongly $\aleph_0$-separated AD family which is not almost-normal under \textsf{CH}, answering a question from Oliveira-Rodrigues and Santos-Ronchim \cite{de2020almost}.\\

We will say that an AD family $\cA$ satisfies a topological property $P$, if and only if $\Psi(\cA)$ does. Given a set $X$ and a cardinal $\kappa$, we denote by $[X]^\kappa$ and $[X]^{<\kappa}$ the set of subsets of $X$ of size $\kappa$ and $<\kappa$, respectively. Also, $[X]^{\leq\kappa}=[X]^\kappa\cup[X]^{<\kappa}$. 
For each $X\subseteq\omega$, $X^0$ will denote $X$ and $X^1$ will denote $\omega\setminus X$. For two infinite subsets $A,B$ of $\omega$, we will say that $A$ meets $B$ if $A\cap B$ is infinite.
We follow \cite{engelking1989general} for topological notation and \cite{kunen2014set} for set theoretic notation. Each undefined weakening of normality can be found in \cite{szeptycki2020weak}.

\section{An almost-normal MAD family}\label{section2}

As we mentioned above, in \cite{szeptycki2020weak}, several counterexamples for the reverse implications in ($\ast$) were given, however, some questions were left open, among them the following two:

\begin{itemize}
    \item Is there an almost-normal not normal AD family?
    \item Is there an almost-normal MAD family?
\end{itemize}

A space $X$ is \emph{almost-normal} (\cite{singal1970almost}) if each pair of closed sets $C,D\subseteq X$, where one of them is regular closed, can be separated. 

Of course, a positive answer for the second question provide a negative answer for the first one. In \cite{de2020almost}, a negative answer for the first question was given. For a subset $X\subseteq 2^\omega$, the AD family $\cA_X\subseteq\cP(2^{<\omega})$ is defined as the family of all sets of the form $\{f\rest n:n\in\omega\}$ with $f\in X$. The result in \cite{de2020almost} is obtained by defining a special class of subsets of $2^\omega$, called almost $Q$-sets, such that $\cA_X$ is the desired family whenever $X$ is an almost $Q$-set and then forcing the set $X$. This result cannot be improved to get MAD since AD families of the form $\cA_X$ are never MAD. 

In this section we will prove that under \textsf{CH} there is an almost-normal MAD family, consistently answering the second question above and improving a result in \cite{alshammari2020quasi}, where the authors prove that there is a quasi-normal MAD family under the same assumption (this result was previously improved in \cite{szeptycki2020weak}, adding the property that the family is Luzin).

\begin{dfn}
Let $\cA$ be an almost disjoint family. A set $D\in[\omega]^\omega$ is a \emph{partitioner} for $\cA$, if for every $A\in\cA$ either, $A\subseteq^*D$ or $A\cap D$ is finite.\\
We will say that two disjoint subfamilies $\cB,\cC\subseteq\cA$ \emph{can be separated}, if there is a partinioner $D$ for $\cA$, such that $B\subseteq^* D$ for every $B\in\cB$ and $|C\cap D|<\omega$ for every $C\in\cC$. In this case, we will say that $D$ is a partitioner for $\cB$ and $\cC$.
\end{dfn}

Notice that if $D$ is a partitioner, $\omega\setminus D$ is a partitioner as well, where the properties of ``almost contained'' and ``is almost disjoint'' have been exchanged. Thus, we can always decide which part of our family is almost contained in the partinioner.

It is known that an AD family $\cA$ is normal, if and only if for each $\cB\subseteq\cA$, $\cB$ and $\cA\setminus\cB$ can be separated \cite{tall1977set}. The respective result for almost normality also holds. We will need the following easy observation.

\begin{fact}
Let $\cA$ be an AD family. For every regular closed set $K\subseteq \Psi(\cA)$, $K=C\cup\{A\in\cA:|A\cap C|=\omega\}$ with $C=K\cap\omega$. 
\end{fact}

\begin{propo}\label{almostnormalequivalence}
An AD family $\cA$ is almost-normal iff for every $C\in[\omega]^\omega$, there exists a partitioner for $\cB=\{A\in\cA:|A\cap C|=\omega\}$ and $\cA\setminus\cB$.
\end{propo}

\begin{proof}
Assume $\cA$ is almost-normal and let $C\subseteq\omega$. Let $$\cB=\{A\in\cA:|A\cap C|=\omega\}.$$
Then $K=\cB\cup C$ is regular closed and $\cA\setminus\cB$ is closed in $\Psi(\cA)$. Since $\cA$ is almost-normal, we can find open disjoint subsets $U,V\subseteq\Psi(\cA)$ such that $K\subseteq U$ and $\cA\setminus\cB\subseteq V$. Define $D=U\cap\omega$ and let $B\in\cB$. Since $U$ contains a basic neighborhood of $B$, it follows that $B\subseteq^* D$.
On the other hand, if $A\in\cA\setminus\cB$, there exists a basic neighborhood of $A$ contained in $V$, thus $A\subseteq^* V\cap\omega$ and therefore $|A\cap D|<\omega$.

Now suppose that each pair $\cB$, $\cA\setminus\cB$ as in the proposition can be separated. Let $F,K\subseteq\Psi(\cA)$ be two disjoint closed sets with $K$ regular closed. There exist $C\subseteq\omega$ such that $K=C\cup\cB$ with $\cB=\{A\in\cA:|A\cap C|=\omega\}$. Let $D$ be a partitioner for $\cB$ and $\cA\setminus\cB$ (where the elements of $\cB$ are those which are almost contained in $D$) and let $E=F\cap\omega$. Define $U=(\cB\cup C\cup D)\setminus E$.

\textbf{Claim:} U is clopen.\\
Given that $\omega$ is discrete, we only care about the points in $\cA$. Let $A\in\cA\setminus\cB$. Since $A\notin\cB$ and $|A\cap D|<\omega$, it follows that $\{A\}\cup(A\setminus(C\cup D))$ is a basic neighborhood of $A$ disjoint from $U$. Then $U$ is closed. If $B\in\mathcal\cB$, $|B\cap E|<\omega$ (otherwise $B\in F$) and $B\subseteq^* D$. Then $\{B\}\cup(D\setminus E)\subseteq U$ contains a basic neighborhood of $B$ showing that $U$ is open.

Finally note that $U$ is a clopen subset disjoint from $F$ and $K\subseteq U$. Thus $\cA$ is almost normal.
\end{proof}

A very related notion on AD families called weakly separation was considered in \cite{brendle1999dow} and \cite{dow1997compact}. Given $\cB,\cC\subseteq\cA$, we say that $D\in[\omega]^\omega$ \emph{weakly separates} $\cB$ and $\cC$, if $D$ meets $B$ for every $B\in\cB$ and $D\cap C$ is finite for every $C\in\cC$. An AD family is \emph{weakly separated} if for any two disjoint subfamilies $\cB,\cC\subseteq\cA$, there is a set $D\in[\omega]^\omega$ that weakly separates $\cB$ and $\cC$. It follows easily that an AD family is normal iff it is almost-normal and weakly separated.

Before the next definition, let us pointing out that for every finite family of partinioners $\{D_i:i<n\}$ of an AD family $\cA$, each boolean combination $\bigcap_{i<n}D_i^{f(i)}$ coded by a function $f:n\to 2$ is also a partinioner for $\cA$.

\begin{dfn}
Let $\cA$ be an almost disjoint family and $\cD\subseteq[\omega]^\omega$ be a family of partitioners for $\cA$. We will say that $\cD$ is a \emph{nice family of partitioners}, if for every $\{D_i:i<n\}\in[\cD]^{<\omega}$ and every $f:n\to2$
$$|\{A\in\cA:A\subseteq^* \bigcap_{i<n}D_i^{f(i)}\}|<\omega\Rightarrow$$
$$\bigcap_{i<n}D_i^{f(i)}=^*\bigcup\lbrace A\in\cA:A\subseteq^* \bigcap_{i<n}D_i^{f(i)}\rbrace$$
\end{dfn}

\begin{lem}\label{nicelemma}
Let $\cA$ be a countable AD family, $C\subseteq\omega$ and let $\cD\subset[\omega]^\omega$ be a countable nice family of partitioners for $\cA$. Then, there exists a partinioner for $\cB=\{A\in\cA:|A\cap C|=\omega\}$ and $\cA\setminus\cB,$ such that $\cD\cup\{D\}$ is a nice family of partitioners for $\cA$.
\end{lem}

\begin{proof}
Enumerate $\cD$ as $\{D_m:m\in\omega\}$ and let $Fn(\omega,2)$ be the set of all finite partial functions $s;\omega\to2$. For each $s\in Fn(\omega,2)$ define
$$D_s=\bigcap_{i\in dom(s)}D_i^{s(i)}.$$
and also define $\cA_s=\{A\in\cA:A\subseteq^*D_s\}$ with a partition of $\cA_s$ into two pieces $\cA_s^+$ and $\cA_s^-$ as follows:
$$\cA_s^+=\{A\in\cA_s:|A\cap C|=\omega\}\textnormal{ and}$$
$$\cA_s^-=\{A\in\cA_s:|A\cap C|<\omega\}.$$

We can enumerate $Fn(\omega,2)$, $\cB=\{A\in\cA:|A\cap C|=\omega\}$ and $\cA\setminus\cB$ as $\{s_n:n\in\omega\}$, $\{B_n:n\in\omega\}$ and $\{B_n':n\in\omega\}$ respectively.

For each $n\in\omega$ define $X_n^+,X_n^-\subseteq\omega$ according to the following cases:

\textit{Case 1: Either, both $\cA_{s_n}^+$ and $\cA_{s_n}^-$ are infinite or both $\cA_{s_n}^+$ and $\cA_{s_n}^-$ are finite.} In this case simply define $X_n^+=\emptyset=X_n^-$.

\textit{Case 2: $\cA_{s_n}^+$ is infinite and $\cA_{s_n}^-$ is finite.} In this case, define $X_n^-=\bigcup\cA_{s_n}^-$ and $X_n^+=D_{s_n}\setminus\bigcup\cA_{s_n}^-$.

\textit{Case 3: $\cA_{s_n}^+$ is finite and $\cA_{s_n}^-$ is infinite.} In this case, define $X_n^+=\bigcup\cA_{s_n}^+$ and $X_n^-=D_{s_n}\setminus\bigcup\cA_{s_n}^+$.

Notice that the family $\{D_{s_n}\setminus \cup\cF:n\in\omega\land\cF\in[\cA]^{<\omega}\}$ is also a nice family of partitioners. Besides, each finite union of elements of $\cA$ is clearly a partitioner and it is easy to see that the family
\begin{itemize}
    \item[($\ast$)]\label{booleanpartitioners}$\{D_{s_n}\setminus \cup\cF:n\in\omega\land\cF\in[\cA]^{<\omega}\}\bigcup\{\cup\cF:\cF\in[\cA]^{<\omega}\}$
\end{itemize}
is still a nice family of partitioners which contains each $X_n^+$ and $X_n^-$. We point out two properties that follow directly from the definitions of the partitioners $X_n^+$ and $X_n^-$ which will be useful later:

\begin{enumerate}[(a)]
    \item\label{X1} $B\cap X_n^-$ is finite for every $B\in\cB$ and $n\in\omega$.
    \item\label{X2} $B'\cap X_n^+$ is finite for every $B'\in\cA\setminus\cB$ and $n\in\omega.$
\end{enumerate}

Define

$$D=\bigcup_{n\in\omega}\left(\left(B_n\cup X_n^+\right)\setminus\bigcup_{i<n}\left(B_i'\cup X_i^-\right)\right).$$

Let $B_n\in\cB$. Thus, by (\ref{X1}) and since $|B_n\cap B'_m|<\omega$ for every $m\in\omega$, 
$$B_n\subseteq^*(B_n\cup X_n^+)\setminus\bigcup_{i<n}\left(B_i'\cup X_i^-\right)\subseteq D.$$

Similarly using (\ref{X2}), we have that $|B_n'\cap D|<\omega$ for every $n\in\omega$. So, $D$ is a partitioner for $\cB$ and $\cA\setminus\cB$.

To finish the proof, we will show that $\cD\cup\{D\}$ is a nice family of partitioners. Given $s\in Fn(\omega,2)$ and $i\in2$ such that $|\{A\in\cA:A\subseteq^* D_{s}\cap D^i\}|<\omega$, if $\cA_s$ is already finite, we have that $D_{s}\cap D^i=^*\{A\in\cA:A\subseteq^* D_{s}\cap D^i\}$ since $D$ is a partitioner. Moreover, $D$ is a partitioner for $\cB$ and $\cA\setminus\cB$ and this implies that either $\{A\in\cA:A\subseteq^* D_{s}\cap D^i\}=\cA_s^+$ if $i=0$ or $\{A\in\cA:A\subseteq^* D_{s}\cap D^i\}=\cA_s^-$ if $i=1$.

Then, we can assume that $\cA_s$ is infinite, and either, $\cA_s^+$ is finite if $i=0$ or $\cA_s^-$ is finite if $i=1$. Let $n\in\omega$ such that $s=s_n$. Assume $i=0$, hence we defined $X_n^+$ and $X_n^-$ according to case 3 and $X_n^+=\bigcup\cA_{s_n}^+$ is a finite union of elements of $\cB$. Thus $X_n^+$ is almost disjoint from  $X_m^-$ and $B_m'$ for every $m\in\omega$. It follows that
$$\bigcup\left(\cA_s^+\right)\subseteq^*(B_n\cup X_n^+)\setminus\bigcup_{i<n}\left(B_i'\cup X_i^-\right)\subseteq D$$
and clearly $\bigcup\left(\cA_s^+\right)\subseteq^* D_s$. So that
$$\bigcup\left(\cA_s^+\right)\subseteq^*D\cap D_s.$$
By case 3, $X_n^-=D_{s_n}\setminus\bigcup\cA_{s_n}^+$. Thus
$$D_s\cap D\setminus\left(\cup\cA_s^+\right)\subseteq\bigcup_{i<n}\left(B_n\cup X_n^+\right).$$

Suppose that $D_s\cap D\setminus\left(\cup\cA_s^+\right)$ is infinite and let $Y\in[D_s\cap D\setminus\left(\cup\cA_s^+\right)]^\omega$. Without loss of generality we can assume that either 
$$Y\in\left[\left( D_s\cap D\setminus\left(\cup\cA_s^+\right)\right)\cap B_i\right]^\omega$$
or 
$$Y\in\left[\left( D_s\cap D\setminus\left(\cup\cA_s^+\right)\right)\cap X_i^+\right]^\omega$$
for some $i<n$. In the first case we would have $Y\in[D_s\cap B_i]^\omega$ which implies that $B_i\subseteq^* D_s$ and then $B_i\in\cA_s^+$, a contradiction. Then we have that $Y\in[D_s\cap D\setminus\left(\cup\cA_s^+\right)\cap X_i^+]^\omega$ for some $i<n$. If $X_i^+$ was defined by case 3, it is a finite union of elements of $\cB$, let's say that $X_i^+=\bigcup_{j\in\sigma}B_j$ for some $\sigma\in[\omega]^{<\omega}$. Then we can repeat the previous argument with some $B_j$.

Thus we can assume that $X_i^+$ was defined by case 2, but even in this situation there is a simple subcase. If $D_s\cap X_i^+$ almost contains only finitely many elements of $\cB$ (in consequence also of $\cA$), it follows by ($\ast$) that $D_s\cap X_i^+$ is a finite union of elements of $\cB$ and we can repeat the previous argument with some $B_j\subseteq D_s\cap X_i^+$.

We then deal with the remaining case and $D_s\cap X_i^+$ almost contains infinitely many elements of $\cB$. In particular, $D_s$ almost contains infinitely many elements of $\cB$ and thus $\cA_s^+$ is infinite, a contradiction. Therefore 
$$D_s\cap D=^*\left(\bigcup\cA_s^+\right)$$
The case $i=1$ is dual and so $\cD\cup\{D\}$ is a nice family of partitioners.
\end{proof}

Recall that given ad AD family $\cA$, the ideal generated by $\cA\cup\omega$ is denoted by $\cI(\cA)$, and $\cI^+(\cA)=\cP(\omega)\setminus\cI(\cA)$ consists of those subsets of $\omega$ that can not be coveret by a finite union of elements of $\cA$ and a finite set. 

\begin{lem}\label{Aconstruction}
Let $\cA$ be a countable AD family, $X\in[\omega]^\omega$ almost disjoint with $\cA$ and let $\cD=\{D_n:n\in\omega\}$ be a nice family of partitioners for $\cA$. Assume that for each $n\in\omega$ there exists $C_n\in[\omega]^\omega$ such that  $D_n$ is a partitioner for $\cB=\{A\in\cA:|A\cap C_n|=\omega\}$ and $\cA\setminus\cB$. Then there exists $A\in[\omega]^\omega$ such that $|A\cap X|=\omega$, $\cA\cup\{A\}$ is AD, each $D_n$ is a partitioner for $\cA\cup\{A\}$ and $A\subseteq^*D_n$ iff $|A\cap C_n|=\omega$ for every $n\in\omega$.
\end{lem}

\begin{proof}
For every $s\in 2^{<\omega}$ let 
$$D_s=\bigcap_{i<|s|}D_{i}^{s(i)}.$$
We can recursively construct $f\in2^\omega$ such that $|D_{f\rest n}\cap X|=\omega$ for every $n\in\omega$. We have that $\{A\in\cA:A\subseteq^*D_{f\rest n}\}$ is infinite for every $n\in\omega$. Otherwise, since $\cD$ is a nice family of partitioners, $X\cap D_{f\rest n}=^*X\cap \bigcup_{i<k} A_i$ for some finite subfamily $\{A_i:i<k\}\subseteq\cA$, which contradicts that $X$ is AD with $\cA$. Thus, $X\cap D_{f\rest n}\in \cI^+(\cA)$. The same argument shows that $C_i\cap D_{f\rest n}\in\cI^+(\cA)$ for every $i<n$ with $f(i)=0$, since each $A\in\cA$ such that $A\subseteq^* D_{f\rest n}$ meets $C_i$ whenever $f(i)=0$. Hence, we can define 
$$F_n\in\left[D_{f\rest n}\setminus\bigcup_{i<n}A_i\right]^{<\omega}$$
such that $F_n\cap X\neq\emptyset$ and $F_n\cap C_i\neq\emptyset$ for every $i<n$ with $f(i)=0$. Even more, we can choose $F_n$ in such a way that $F_n\cap F_m=\emptyset$ for every $n\neq m$. 

Define $A=\bigcup_{n\in\omega} F_n$. It is clear that $|A\cap X|=\omega$ and $\cA\cup\{A\}$ is AD from the definition of the $F_n$'s. Now assume that $f(n)=1$, then $F_k\subseteq D_{f\rest k}\subseteq \omega\setminus D_n$ for every $k>n$ which implies that $|A\cap D_n|<\omega$. Similarly $A\subseteq^* D_n$ whenever $f(n)=0$, and in this case, $F_k\cap C_n\neq\emptyset$ for every $k>n$. Therefore $|A\cap C_n|=\omega$. This finishes the proof.
\end{proof}

\begin{teo}\label{Teorema}{\textsf{(CH)}}
There is an almost-normal MAD family.
\end{teo}

\begin{proof}
Enumerate $[\omega]^\omega=\{X_\alpha:\alpha<\omega_1\}$ with $X_n=\omega$ for every $n\in\omega$. We will recursively construct a MAD family $\cA=\{A_\alpha:\alpha<\omega_1\}$ and a family of partitioners $\cD=\{D_\alpha:\alpha<\omega_1\}$ such that if $\cA_\alpha=\{A_\beta:\beta<\alpha\}$ and $\cD_\alpha=\{D_\beta:\beta<\alpha\}$ then:
\begin{enumerate}
    \item $\cA_\alpha$ is AD.
    \item\label{p1} If $X_\alpha$ is AD with $\cA_\alpha$, then $|A_\alpha\cap X_\alpha|=\omega$.
    
    \item $D_\alpha$ is a partitioner for $\cB=\{A\in\cA_\alpha:|A\cap X_\alpha|=\omega\}$ and $\cA_\alpha\setminus\cB$.
    
    \item\label{p3} Either, $A_\alpha\subseteq^* D_\beta$ or $|A_\alpha\cap D_\beta|<\omega$ for every $\beta\leq\alpha$.
    
    \item\label{p4} $A_\alpha\subseteq^*D_\beta$ iff $|A_\alpha\cap X_\beta|=\omega$.
    
    \item $\{D_\beta:\beta\leq\alpha\}$ is a nice family of partitioners for $\cA_\alpha$.
\end{enumerate}

Let $\{A_n:n\in\omega\}\subseteq[\omega]^\omega$ be a partition of $\omega$ into infinite pieces and define $D_n=\omega$ for every $n\in\omega$. This family clearly satisfies the above conditions. Assume we have constructed $\cA_\alpha$ and $\cD_\alpha$ as above. We can apply lemma \ref{nicelemma} to the triple $(\cA_\alpha,\cD_\alpha,X_\alpha)$ to obtain $D_\alpha$.

For the construction of $A_\alpha$, let $X=X_\alpha$ if $X_\alpha$ is AD with $\cA_\alpha$, otherwise let $X$ be any infinite subset of $\omega$ almost disjoint with $\cA_\alpha$. By lemma \ref{Aconstruction} applied to $\cA_\alpha$, $\{D_\beta:\beta\leq\alpha\}$ (with their respective $C_\beta=X_\beta$) and $X$, we can find $A_\alpha$ as required.

It is clear from point (\ref{p1}) that $\cA$ is MAD. Also, if $C\in[\omega]^\omega$, there exists $\alpha<\omega_1$ such that $C=X_\alpha$. Hence $D_\alpha$ is a partitioner for $\cA_\alpha$ as in proposition \ref{almostnormalequivalence}. Moreover, point (\ref{p3}) and point (\ref{p4}) ensure that $D_\alpha$ is preserved for $\beta\geq\alpha$. Thus $D_\alpha$ is a partitioner for $\cA$ as required in proposition \ref{almostnormalequivalence} with $C=X_\alpha$. We can conclude that $\cA$ is almost-normal.
\end{proof}

It was mentioned before that in \cite{szeptycki2020weak}, a quasi-normal Luzin MAD family was constructed, then it is natural to ask the following question:

\begin{que}
(\textsf{}CH) Is there a Luzin MAD family which is almost-normal?
\end{que}

\section{There may be no almost-normal MAD families}\label{section3}

There are many reasons for which one could think that it is not possible to obtain Theorem \ref{Teorema} without assuming \textsf{CH}. The most obvious reason is that after $\omega_1$-many steps, we have already constructed a Luzin family $\cA$. Then, we can not get a partitioner as in Proposition \ref{almostnormalequivalence} for a given set $C\subseteq[\omega]^\omega$, whenever it meets uncountable many elements of $\cA$ and it is almost disjoint from uncountable many elements of $\cA$ as well. Indeed, this situation could be unavoidable as we will see below.

Recall that the Proper Forcing Axiom (\textsf{PFA}) is the assertion that for every proper forcing $\mathbb{P}$ and every family $\cD$ of $\omega_1$-many open dense subsets of $\mathbb{P}$ there exists a $\cD$-generic filter for $\mathbb{P}$. If we replace ``proper'' by ``ccc'' and ``$\omega_1$'' by ``$<\mathfrak{c}$'' we get the definition of Martin's Axiom (\textsf{MA}). It is well known that \textsf{PFA} implies \textsf{MA}$+\mathfrak{c}=\omega_2$. Under \textsf{PFA} we can not avoid the existence of Luzin subfamilies due to the following result.

\begin{teo}\cite{dow2011sequential}\label{pfa}
Each MAD family contains a Luzin subfamily.
\end{teo}

The existence of a set $C$ as above, that ``wants to separate'' the Luzin subfamily is also insured by the next theorem.

\begin{teo}\cite{martin1970internal}\label{Cohen}
(\textsf{MA}) For every pair of families $\cA,\cB\subseteq[\omega]^\omega$ of size $<\mathfrak{c}$ such that for every $K\in[\cA]^{<\omega}$ and $B\in\cB$, $B\setminus\bigcup K$ is infinite, there exists $C\in[\omega]^\omega$ such that $C\cap A$ is finite for every $A\in\cA$ and $C$ meets $B$ for every $B\in\cB$.
\end{teo}

Now it follows easily that there are no almost-normal MAD families in the presence of \textsf{PFA}.

\begin{cor}\label{coro}
\textsf{PFA} implies that there are no almost-normal MAD families.  
\end{cor}

\begin{proof}
Let $\cA$ be a MAD family and let $\cA'\subseteq\cA$ be a Luzin subfamily. We can split $\cA'$ into two uncountable disjoint subfamilies $\cB,\cC\subseteq\cA'$. By Theorem \ref{Cohen}, we can find a set $X\subseteq\omega$ that weakly separates $\cB$ and $\cC$, that is, $X\cap C$ is finite for every $C\in\cC$ and $X$ meets $B$ for every $B\in\cB$. Thus, $\cK=\{A\in\cA:|X\cap A|=\omega\}$ and $\cA\setminus\cK$ cannot be separated since $\cB\subseteq\cK$, $\cC\subseteq\cA\setminus\cK$ and $\cB,\cC$ are uncountable subfamilies of a Luzin family. Therefore $\cA$ is not almost-normal.
\end{proof}

In \cite{dow2012martin}, it is shown that it is consistent with \textsf{MA} that there is a MAD family which contains no Luzin subfamilies, it could be possible that the only thing that blocks the existence of almost-normal MAD families is the existence of Luzin subfamilies, so we ask the following:

\begin{que}
Is it consistent with \textsf{MA} that there are almost-normal MAD families?
\end{que}

\section{Partly-normal not quasi-normal AD families}\label{section4}

In this section, we will consider the next question stated in \cite{szeptycki2020weak} and will provide a positive answer.
\begin{itemize}
    \item Is there a partly-normal not quasi-normal AD family?
\end{itemize}

We will say that a space $X$ is \emph{partly-normal} if any pair of disjoint closed sets $A,B\subseteq X$, where $A$ is regular closed and $B$ is $\pi$-closed (a finite intersection of regular closed sets), can be separated \cite{alshammari2020quasi}. A space $X$ is \emph{quasi-normal} if any two disjoint $\pi$-closed sets can be separated \cite{zaitsev1968some}.

Most of the examples in \cite{szeptycki2020weak} were constructed using AD families of true cardinality $\mathfrak{c}$. 
For an AD family $\cA$ and $W\subseteq\omega$, we will denote by $\cA\rest W$ the set of $A\in\cA$ such that $A$ meets $W$. 
An AD family is of \emph{true cardinality $\mathfrak{c}$}, if for every $W\subseteq\omega$, either, $\cA\rest W$ is finite or has size $\mathfrak{c}$. 
It is well known that the existence of (M)AD families of true cardinality $\mathfrak{c}$ is equivalent to the existence of completely separable (M)AD families. An AD family is \emph{completely separable}, if for any $X\subseteq\omega$ such that $|\cA\rest X|=\omega$, there is an $A\in\cA$ such that $A\subseteq X$ \cite{hechler1971classifying}. While completely separable AD families do exist in \textsf{ZFC}, the existence of completely separable MAD families in \textsf{ZFC}, asked first by Erd\"os and Shelah \cite{erdos1972separability}, is one of the more interesting and central questions concerning almost disjoint families. Completely separable MAD families were constructed under several assumptions, $\mathfrak{a}=\mathfrak{c},\mathfrak{b}=\mathfrak{d},\mathfrak{d}\leq\mathfrak{a}$ and $\mathfrak{s}=\omega_1$ (see \cite{balcar1984almost,balcar1989disjoint,simon1996note}). 
Then Shelah showed that they exists under $\mathfrak{s}<\mathfrak{a}$, and also under $\mathfrak{s}\geq\mathfrak{a}$ assuming covering type assumptions \cite{shelah2011mad}, in particular, they exists under $\mathfrak{c}<\aleph_\omega$ (see \cite{hruvsak2014almost}). Later work of Mildenberger, Raghavan and Steprans \cite{mildenberger2014splitting}, showed that the covering type assumption is not needed in the case $\mathfrak{s}=\mathfrak{a}$, and then completely separable MAD families exists under $\mathfrak{s}\leq\mathfrak{a}$. 
The existence of AD families of true cardinality $\mathfrak{c}$ is particularly useful for constructions of AD families with strong combinatorial properties, since they usually need recursive constructions of length continuum (see, for example, \cite{raghavan2010almost}). We will use an AD family of true cardinality $\mathfrak{c}$ to construct a partly-normal not quasi-normal AD family. 
First, observe that we can always assume that an infinite AD family $\cA$, contains an infinite partition of $\omega$ into infinite pieces, since we can take $\{A_n:n\in\omega\}\subseteq\cA$ and substitute $A_n$ by $A_n'=(A_n\cup\{n\})\setminus\bigcup_{i<n}A_i'$.

We will use the notion of dominance of function in $\baire$. Recall that for two functions $f,g\in\baire$, we say that $f<^*g$, if the set $\{n\in\omega:f(n)\geq g(n)\}$ is finite. We say that a function $g$ dominates a family $\cF\subseteq\baire$, if $g>^*f$ for every $f\in\cF$. A family $\cB\subseteq\baire$ is an unbounded family, if no single function $f\in\baire$ dominates $\cB$. A family $\cF\subseteq\baire$ is a dominating family, if for every $g\in\baire$, there is an $f\in\cF$ such that $f>^*g$. The unbounding number $\mathfrak{b}$ is the least size of an unbounded family, and the dominating number $\mathfrak{d}$ is the least size of a dominating family. It is well known that $\omega<\mathfrak{b}\leq\mathfrak{d}\leq\mathfrak{c}$ (see \cite{blass}).
We are now ready to prove the following result.

\begin{teo}
There is a partly-normal AD family which is not quasi-normal. 
\end{teo}

\begin{proof}
Partition $\omega$ in four infinite sets $W_0,W_1,V_0$ and $V_1$. Further, partition both, $W_0$ and $W_1$ into infinitely many infinite sets, that is,
$$W_0=\bigcup_{n\in\omega}P_n$$
and
$$W_1=\bigcup_{n\in\omega}Q_n.$$
Let $\cA_{W_0},\cA_{W_1},\cA_{V_0}$ and $\cA_{V_1}$ be AD families of true cardinality $\mathfrak{c}$ in each of the four sets $W_0,W_1,V_0$ and $V_1$. We can assume that $\{P_n:n\in\omega\}\subseteq\cA_{W_0}$ and $\{Q_n:n\in\omega\}\subseteq\cA_{W_1}$. 
We will recursively construct our family putting together some elements of these AD families of true cardinality $\mathfrak{c}$. For ease of notation, let $\cE=\cA_{W_0}\cup\cA_{W_1}\cup\cA_{V_0}\cup\cA_{V_1}$.

For every $n\in\omega$, define $A_n=P_n\cup Q_n$. Let $\{f_\alpha:\alpha<\mathfrak{c}\}$ be a dominating family of functions in $\omega^\omega$. List all pairs $(C,\cD)$ where $\cD\in\left[\left[\omega\right]^\omega\right]^{<\omega}$ and $C\in[\omega]^\omega$ as $\{(C_\alpha,\cD_\alpha):\alpha<\mathfrak{c}\}$. We will built finite sets $\cF_\alpha^0,\cF_\alpha^1\in[\cE]^{<\omega}$ recursively, so that each $\cF_\alpha^0$ will contain exactly one element of each family  $\cA_{W0}, \cA_{V_0}$ and $\cA_{V_1}$, and each $\cF_\alpha^1$ will intersects exactly two of the families $\cA_{W_0},\cA_{W_1},\cA_{V_0},\cA_{V_1}$. In particular, no element of the form $A_n$ or $\cup\cF_\alpha^i$ will meet the four sets $W_0,W_1,V_0$ and $V_1$

Assume we have constructed $\cF_\beta^0$ and $\cF_\beta^1$ for $\beta<\alpha$. For $\cF_\alpha^0$, consider the set $X=\bigcup_{n\in\omega}(P_n\setminus f_\alpha(n))$. Since $X$ meets infinitely many elements of $\cA_{W_0}$, it follows that $X$ meets $\mathfrak{c}$-many elements of $\cA_{W_0}$. Choose $$A\in\cA_{W_0}\setminus\left(\bigcup_{\beta<\alpha}(\cF_\beta^0\cup\cF_\beta^1)\cup\{P_n:n\in\omega\}\}\right)$$
such that $A$ meets $X$. Also pick $B\in\cA_{V_0}\setminus\bigcup_{\beta<\alpha}(\cF_\beta^0\cup\cF_\beta^1)$ and $C\in\cA_{V_1}\setminus\bigcup_{\beta<\alpha}(\cF_\beta^0\cup\cF_\beta^1)$ arbitrary and define $\cF_\alpha^0=\{A,B,C\}$.

For $\cF_\alpha^1$ consider the pair $(C_\alpha,\cD_\alpha)$ and let $\cD_\alpha=\{D_\alpha^j:j<n\}$. Define 
$$\cC=\{A\in\cE:A\textnormal{ meets }C_\alpha\}$$
and 
$$\cB=\{A\in\cE:\forall\ j<n\ (A\textnormal{ meets }D_\alpha^j)\}.$$

If either $\cB$ or $\cC$ are finite, simply define $\cF_\alpha^1=\emptyset$. Otherwise, we have some cases. Since $\cB$ and $\cC$ are infinite, there are $Y,Z\in\{W_0,W_1,V_0,V_1\}$ such that $\cB\rest Y$ and $\cC\rest Z$ are infinite. Notice that since $\cB\rest Y$ is infinite, in particular, there are infinitely many elements of $\cA_Y$ which meet $D_\alpha^j$ for every $j<n$. Thus, there are $\mathfrak{c}$-many of these elements. Pick, for every $j<n$, an element $B_j\in\cA_Y$ such that $B_j$ meets $D_\alpha^j$ and $$B_j\in\cA_Y\setminus\left(\cF_\alpha^0\cup\bigcup_{\beta<\alpha}\left(\cF_\beta^0\cup\cF_\beta^1\right)\right).$$
Notice that $\bigcup_{j<n}B_j\subseteq Y$. Similarly, since $\cC$ is infinite, we can find a set $C'\in\cA_X$ such that $C'$ meets $C_\alpha$ and
$$C'\in\cA_X\setminus\left(\cF_\alpha^0\cup\bigcup_{\beta<\alpha}\left(\cF_\beta^0\cup\cF_\beta^1\right)\right).$$
Define $\cF_\alpha^1=\{C'\}\cup\{B_j:j<n\}$. This finishes the construction of the $\cF_\alpha^i$'s.

Now, we can describe our AD family. Let 
$$\cA=\{A_n:n\in\omega\}\bigcup\{\cup\cF_\alpha^i:\alpha<\mathfrak{c}\ \land\  i<2\}.$$
It is clear that it is AD since each of its elements is a finite union of elements of $\cE$.

To see that it is partly-normal, let $K_0$ and $K_1$ be two disjoint closed subsets of $\Psi(\cA)$ such that $K_0$ is regular closed and $K_1$ is $\pi$-closed. Hence, there are $C\in[\omega]^\omega$ and $\{D_j:j<n\}\subseteq[\omega]^\omega$ such that $K_0=\overline{C}$ and $K_1=\bigcap_{j<n}\overline{D_j}$ in $\Psi(\cA)$. Suppose that $K_0\cap\cA$ is finite, thus
$$U=C\cup\bigcup\lbrace\lbrace A\rbrace\cup(A\setminus\cap_{j<n}D_j):A\textnormal{ meets }C\rbrace$$
is a clopen subset which separates $K_0$ and $K_1$. A similar argument shows that if $K_1\cap\cA$ is finite, we can separate $K_0$ and $K_1$. 

We can then assume that both $K_0\cap\cA$ and $K_1\cap\cA$ are infinite. Let $\alpha<\mathfrak{c}$ such that $(C,\{D_j:j<n\})=(C_\alpha,\cD_\alpha)$. By the previous assumption, $\cF_\alpha^1$ is not empty. So, $\cF_\alpha^1=\{C'\}\cup\{B_j:j<n\}$, where $C'$ meets $C$ and $B_j$ meets $D_j$ for every $j<n$, which implies that $\cup\cF_\alpha^1\in K_0\cap K_1$, a contradiction. Hence the case where $K_0\cap\cA$ and $K_1\cap\cA$ are infinite is not possible. 

To see that $\cA$ is not quasi-normal, consider $W=\overline{W_0}\cap\overline{W_1}$ and $V=\overline{V_0}\cap\overline{V_1}$. These two closed sets are disjoint since no element of $\cA$ intersects the four sets $W_0,W_1,V_0,V_1$, which are a partition of $\omega$. Let $U$ be an open set containing $W$. We have that $A_n$ meets both $W_0$ and $W_1$ for every $n\in\omega$, then each $A_n\in W$. In particular, $P_n\subseteq A_n\subseteq^*U$. Let $f\in\omega^\omega$ such that $P_n\setminus f(n)\subseteq U$. We can find an $\alpha<\mathfrak{c}$ such that $f<^*f_\alpha$. Then, at step $\alpha$, we defined $\cF_\alpha^0=\{A,B,C\}$ in such a way that $A$ meets $\bigcup_{n\in\omega}(P_n\setminus f_\alpha(n))$ (and consequently $\bigcup_{n\in\omega}(P_n\setminus f(n))\subseteq U$), $B$ meets $V_0$ and $C$ meets $V_1$. Therefore $\cup\cF_\alpha^0\in K_1$ but no open set containing it can be disjoint from $U$, which makes us impossible to separate $K_0$ from $K_1$.
\end{proof}

All known counterexamples of the normality-like properties considered here, with exception of an almost-normal not normal AD family, can be constructed in \textsf{ZFC} alone. Hence, it is natural to ask if such a space can also exists in \textsf{ZFC}. We already know that no counterexample can be MAD by Corollary \ref{coro}. In \cite{de2020almost}, the cardinal $\mathfrak{an}$ is defined as the least cardinality of an almost-normal not normal AD family, and it is noted that $\mathfrak{ap}\leq\mathfrak{an}$, whenever $\mathfrak{an}$ is well defined, i.e., whenever there is an almost-normal not normal AD family. Here $\mathfrak{ap}$ is defined as the least cardinality of an AD family which is not weakly separated \cite{brendle1999dow}. Since it is consistent that $\mathfrak{ap}=\mathfrak{c}$, the unresolved portion of the question of whether there are almost-normal MAD families, can be stated as follows: 

\begin{que}
Does there exist (in \textsf{ZFC}) an almost-normal AD family which is not normal? (an almost-normal AD family of size $\mathfrak{c}$?)
\end{que}

On the other hand, it was proved in \cite{de2020almost}, that there is an almost-normal not normal AD family of size $\omega_1<\mathfrak{c}$. Hence, even though the first part of the above question might have a ``yes'' answer, the proof may go by cases (in some models all such families have size $<\mathfrak{c}$ while in others, all such families have size $\mathfrak{c}$) and then the second part of the question could have a ``no'' answer.\\

In \cite{hruvsak2013n}, a study on the relation between normality and the existence of Luzin-type subfamilies was done. We call a pair $\cB=\{B_\alpha:\alpha<\omega_1\}$ and $\cC=\{C_\alpha:\alpha<\omega_1\}$ of subfamilies of $[\omega]^\omega$ a \emph{Luzin gap} if there is an $m\in\omega$ such that:
\begin{enumerate}
    \item $A_\alpha\cap B_\alpha\subseteq m$ for every $\alpha<\omega_1$ and
    \item $(A_\alpha\cap B_\beta)\cup (A_\beta\cap B_\alpha)\nsubseteq m$ but $A_\alpha\cap B_\beta$ is finite for every $\alpha\neq\beta<\omega_1.$
\end{enumerate}

It is known that every Luzin family contains many Luzin gaps and if $\cB$ and $\cC$ forms a Luzin gap, they can not be separated. Thus, AD families which contain Luzin gaps are not normal. Moreover, Luzin gaps are indesctructible by forcing notions which preserves $\omega_1$, thus, Luzin gaps can not be normal in any of these forcing extensions. A generaization of Luzin gaps is the following:

\begin{dfn}\cite{hruvsak2013n}
Let $n\in\omega$ and $B_i=\{B_\alpha^i:\alpha<\omega_1\}$ be disjoint subfamilies of an AD family $\cA$ for $i<n$. We say that $\langle B_i:i<n\rangle$ forms an $n$-Luzin gap if there is an $m\in\omega$ such that:
\begin{enumerate}
    \item $B_\alpha^i\cap B_\alpha^j\subseteq m$ for all $i\neq j$, $\alpha<\omega_1$ and
    \item $\bigcup_{i\neq j}(B_\alpha^i\cap B_\beta^j)\nsubseteq m$ for every $\alpha\neq\beta<\omega_1.$
\end{enumerate}
\end{dfn}

Let $P$ be any property of AD families. An AD family is said to be \emph{potentially $P$} \cite{hruvsak2013n}, if there is a forcing notion $\mathbb{P}$, such that $\forces_{\mathbb{P}}$``$\cA$ is $P$''. Hence, an AD family fails to be potentially normal if it contains Luzin gaps. An interesting result arises when $n$-Luzin gaps are considered under \textsf{MA}.

\begin{teo}\cite{hruvsak2013n}
Assume \textsf{MA} and let $\cA$ be an AD family. Then $\cA$ is normal if and only if $|\cA|<\mathfrak{c}$ and $\cA$ does not contain $n$-Luzin gaps for any $n\in\omega$.
\end{teo}

A result in \textsf{ZFC} that could be useful to the study on normality-like properties is the following:

\begin{teo}\cite{hruvsak2013n}
The following are equivalent for an AD family $\cA$:
\begin{enumerate}
    \item $\cA$ does not contain $n$-Luzin gaps for any $n\in\omega$,
    \item $\cA$ is potentially normal,
    \item $\cA$ is potentially $\mathbb{R}$-embeddable.
\end{enumerate}
\end{teo}

Recall that $\cA$ is $\mathbb{R}$-embeddable if there is an injective and continuous function $\varphi:\Psi(\cA)\to\mathbb{R}$. Hence, one could ask the relation between these concepts and the weakenings of normality.

\begin{que}
Are almost-normal AD families potentially normal?
\end{que}

Since it is consistent that there are quasi-normal AD families which contain Luzin families, we can not ask the above question for weaker normality-like properties in \textsf{ZFC}.

\begin{que}
Is it consistent that quasi-normal (partly-normal, mildly-normal) AD families are potentially normal?
\end{que}

\section{On strongly $\aleph_0$-separated AD families}\label{section5}

The concept of strongly $\aleph_0$-separated AD families was introduced in \cite{szeptycki2020weak} by the authors. An AD family $\cA$ is \emph{strongly $\aleph_0$-separated}, if for every two disjoint countable subfamilies $\cB,\cC\subseteq\cA$, there is a partitioner $D$ for $\cA$, such that $B\subseteq^*D$ for every $B\in\cB$ and $C\cap D$ is finite for every $C\in\cC$. There, it was shown that almost-normal AD families  are strongly $\aleph_0$-separated and that there is a strongly $\aleph_0$-separated MAD family under \textsf{CH}.

The requirement of one of the subfamilies being countable was modified in \cite{de2020almost} in order to define a stronger concept: An AD family is \emph{strongly $(\aleph_0,<\mathfrak{c})$-separated}, if for every two disjoint subfamilies $\cB,\cC\subseteq\cA$, where $\cB$ is countable and $|\cC|<\mathfrak{c}$, there is a partitioner $D$ for $\cA$ such that $B\subseteq^*D$ for every $B\in\cB$ and $C\cap D$ is finite for every $C\in\cC$. The relation of these two concepts and almost-normality was studied in \cite{szeptycki2020weak} and \cite{de2020almost}, however, the next question remained unanswered \cite{de2020almost}:

\begin{itemize}
    \item Does \textsf{CH} imply that strongly $\aleph_0$-separated AD families are almost-normal?
\end{itemize}
We will answer this question in the negative. For this purpose, recall that for $S,X\in[\omega]^\omega$, we say that $S$ \emph{splits} $X$, if both $S\cap X$ and $X\setminus S$ are infinite. A family $\cS\subseteq\cP(\omega)$ is a \emph{splitting family}, if for every $X\in[\omega]^\omega$, there exists $S\in\cS$ such that $S$ splits $X$, and $\mathfrak{s}$ is the least size of a splitting family. The splitting number $\mathfrak{s}$ is a cardinal invariant of the continuum, hence $\omega<\mathfrak{s}\leq\mathfrak{c}$. In particular, for every countable family $\cH\subseteq[\omega]^\omega$, there exists $X\in[\omega]^\omega$ which is not split by any element of $\cH$, i.e., either, $X\cap H$ is finite or $X\subseteq^* H$ for every $H\in\cH$ (see \cite{blass}).

\begin{teo}{\textsf{(CH)}}
There is a strongly $\aleph_0$-separated AD family which is not almost-normal.
\end{teo}

\begin{proof}
Let $\{d_\alpha:\alpha<\omega_1\}$ be a dominating family of functions and enumerate all pairs $(a,b)\in[\omega_1]^{\leq\omega}\times[\omega_1]^{\leq\omega}$, such that $a\cap b=\emptyset$ as $\{(a_\alpha,b_\alpha)\colon\omega\leq \alpha<\omega_1\}$. We can assume without loss of generality that $a_\alpha\cup b_\alpha\subseteq\alpha$ for every $\alpha<\omega_1$. Partition $\omega=V\cup W$ into two infinite sets and let $\varphi\colon V\to W$ be a bijection. Moreover, partition $V=\bigcup_{n\in\omega}A_n$, into infinitely many infinite sets.

We will recursively construct $A_\alpha$ and $D_\alpha$ for $\omega\leq\alpha<\omega_1$ such that if $\cA_\alpha=\{A_\beta:\beta<\alpha\}$, then the following holds:

\vspace{3mm}
\begin{enumerate}
    
    \item\label{strongly1} $\cA_\alpha$ is an almost disjoint family.
    
    \item\label{strongly2} $|A_\alpha\cap A_n|\leq 1$ for every $\omega\leq\alpha<\omega_1$ and every $n\in\omega$.
    
    \item\label{strongly3} $A_\alpha$ meets $\bigcup_{n\in\omega}\left(A_n\setminus d_\alpha(n)\right)$.
    
    \item\label{strongly4} $D_\alpha$ is a partitioner for $\cA_\alpha$ such that $A_\beta\subseteq^* D_\alpha$ for every $\beta\in a_\alpha$ and $A_\gamma\cap D_\alpha$ is finite for every $\gamma\in b_\alpha$.
    
    \item\label{strongly5} Either, $A_\alpha\subseteq^*D_\beta$ or $A_\alpha\cap D_\beta$ is finite for every $\beta<\alpha$. 
    
\end{enumerate}
\vspace{3mm}

Suppose we have defined $\cA_\alpha=\{A_\beta:\beta<\alpha\}$ and $\cD_\alpha=\{D_\beta:\omega\leq\beta<\alpha\}$ with the above properties. We shall define $D_\alpha$ and $A_\alpha$.

Consider the pair $(a_\alpha,b_\alpha)$. Let $\cB=\{A_\beta:\beta\in a_\alpha\}$ and $\cC=\{A_\beta:\beta<\alpha\ \land\ \beta\notin a_\alpha\}$. Since $\alpha$ is countable we can enumerate both sets as $\cB=\{B_n:n\in\omega\}$ and $\cC=\{C_n:n\in\omega\}$. Define 

$$D_\alpha=\bigcup_{n\in\omega}\left(B_n\setminus\bigcup_{i<n}C_i\right).$$

Since $\cA_\alpha=\cB\cup\cC$ is AD, it is easy to see that $D_\alpha$ is a partitioner for $\cA_\alpha$ and it follows from the definition that satisfies property (\ref{strongly4}).

Now we turn to the construction of $A_\alpha$. For every infinite ordinal $\beta<\alpha$, there is a function $f_\beta$ such that $A_\beta\cap A_n\subseteq f_\beta(n)$ for every $n\in\omega$. Define for every infinite ordinal $\beta\leq\alpha$,  $H_\beta=\{n\in\omega:D_\beta\textnormal{ meets }A_n\}$. Notice that since $D_\beta$ is a partitioner, $A_n\subseteq^* D_\beta$ whenever $D_\alpha$ meets $A_n$. Thus we can also define a function $g_\beta\in\baire$ such that $A_n\setminus g_\beta(n)\subseteq D_\beta$ if $n\in H_\beta$ and $D_\beta\cap A_n\subseteq g_\beta(n)$ otherwise. Let $r\in\baire$ such that $r$ dominates $\{d_\beta:\beta\leq\alpha\}\cup\{f_\beta:\omega\leq\beta<\alpha\}\cup\{g_\beta:\omega\leq\beta\leq\alpha\}$.

Since the family $\{H_\beta:\omega\leq\beta\leq\alpha\}$ is countable, we can also find a set $H\in[\omega]^\omega$ such that for every infinite ordinal $\beta\leq\alpha$ either, $H\cap H_\beta$ is finite or $H\subseteq^* H_\beta$. For every $n\in H$, let $x_n\in A_n\setminus r(n)$. Define $A_\alpha=\{x_n:n\in\omega\}$.

It is clear that $A_\alpha$ satisfies (\ref{strongly2}), given that the family $A_n$ is a partition of $V$. To see that $A_\alpha$ satisfies (\ref{strongly3}), simply note that $r>^*d_\alpha$. Let $\beta<\alpha$ an infinite ordinal and let $k\in\omega$ such that $r(n)>f_\beta(n)$ for every $n>k$. Then 
$$x_n\in A_n\setminus r(n)\subseteq A_n\setminus f_\beta(n)\subseteq A_n\setminus A_\beta$$
for every $n>k$, showing that $A_\alpha\cap A_\beta$ is finite and then property (\ref{strongly1}) holds as well. For property (\ref{strongly5}), consider the set $H_\beta$. If $H\cap H_\beta$ is finite, we can find $k\in\omega$ such that $H\cap H_\beta\subseteq k$ and $r(n)>g_\beta(n)$ for every $n>k$. Hence, for every $n\in H\setminus k$, $n\notin H_\beta$ implies that $D_\beta\cap A_n\subseteq g_\beta(n)$ and 
$$x_n\in A_n\setminus r(n)\subseteq A_n\setminus g_\beta(n)\subseteq A_n\setminus D_\beta,$$
whence $A_\alpha\cap D_\beta$ is finite. On the other hand, if $H\subseteq^* H_\beta$, we can find $k\in\omega$ such that $H\setminus k\subseteq H_\beta$ and $r(n)>g_\beta(n)$ for every $n>k$. Then, for every $n\in H\setminus k\subseteq H_\beta$, $g_\beta(n)$ was defined so that $A_n\setminus g_\beta(n)\subseteq D_\beta$ and 
$$x_n\in A_n\setminus r(n)\subseteq A_n\setminus g_\beta(n)\subseteq D_\beta,$$
proving that $A_\alpha\subseteq^* D_\beta$. This finishes the recursive construction.

We are going to make a last modification to $\cA=\{A_\alpha:\alpha<\omega_1\}$ in order to get the desired family. For every $\alpha<\omega_1$ define $\overline{A_\alpha}$ as follows:
$$\overline{A_\alpha}=\begin{cases}
A_\alpha & \textnormal{ if }\alpha<\omega\\
A_\alpha\cup\varphi[A_\alpha] & \textnormal{ if }\alpha\geq\omega 
\end{cases}$$
Similarly define $\overline{D_\alpha}=D_\alpha\cup\varphi[D_\alpha]$ for $\omega\leq\alpha<\omega_1$. Since $\varphi$ is a bijection between two disjoint sets $V$ and $W$, if $\overline{\cA}=\{\overline{A_\alpha}:\alpha<\omega_1\}$ and $\overline{\cA_\alpha}=\{\overline{A_\beta}:\beta<\alpha\}$, properties (\ref{strongly1})-(\ref{strongly5}) also hold replacing $A_\tau$ by $\overline{A_\tau}$ and $D_\tau$ by $\overline{D_\tau}$.\\

\textbf{Claim:} $\overline{\cA}$ is strongly $\aleph_0$-separated.\\
Let $\cA_0,\cA_1\in[\overline{\cA}]^{\leq\omega}$ be disjoint subfamilies. Define $a=\{\delta:\overline{A_\delta}\in\cA_0\}$ and $b=\{\beta:\overline{A_\beta}\in\cA_1\}$. There exists $\alpha<\omega_1$ such that $(a,b)=(a_\alpha,b_\alpha)$. Thus $\overline{D_\alpha}$ is a partitioner for $\overline{\cA_\alpha}$ and was chosen so that $\overline{D_\alpha}$ separates $\cA_0$ and $\cA_1$ by property (\ref{strongly4}). Moreover, since $\overline{A_\gamma}$ is either, almost disjoint or almost contained in $\overline{\cD_\alpha}$ for every $\gamma\geq\alpha$, $\overline{D_\alpha}$ is indeed, a partitioner for $\overline{\cA}$ which separates $\cA_0$ and $\cA_1$.\\

\textbf{Claim:} $\overline{\cA}$ is not almost-normal.\\
For every $n\in\omega$, $\overline{A_n}\subseteq V$ which is disjoint from $W$. In addition, $\overline{A_\alpha}\cap W=\varphi[A_\alpha]$ is an infinite set. It suffices now to prove that $\{\overline{A_n}:n\in\omega\}$ and $\{\overline{A_\alpha}:\omega\leq\alpha<\omega_1\}$ can not be separated. Let $D$ be such that $\overline{A_n}\subseteq^*D$ for every $n\in\omega$. There exists $f\in\baire$ such that $\overline{A_n}\setminus f(n)\subseteq D$. Choose $\alpha<\omega_1$ such that $d_\alpha>^*f$. Then $\overline{A_\alpha}$ meets $\bigcup_{n\in\omega}\left(\overline{A_n}\setminus d_\alpha(n)\right)\subseteq^*D$. Hence, $D$ is not a partitioner for $\{\overline{A_n}:n\in\omega\}$ and $\{\overline{A_\alpha}:\omega\leq\alpha<\omega_1\}$.
\end{proof}

We have answered Question 7.3 from \cite{de2020almost} in the negative, in particular, under \textsf{CH}, there is a strogly $(\aleph_0,<\mathfrak{c})$-separated AD family which is not almost-normal. This result also follows from \textsf{PFA}, actually, something stronger is true. Let $P$ be a given property. We will say that MAD families with property $P$ \emph{exists generically} if all AD families of size less than $\mathfrak{c}$ can be extended to a MAD family with property $P$. Generic existence of MAD families was introduced in \cite{guzman2017generic} and it was proved in \cite{de2020almost} that under $\mathfrak{b}=\mathfrak{c}=\mathfrak{s}$, completely separable MAD families which are strongly $(\aleph_0,<\mathfrak{c})$-separated exist generically. Since the hypothesis hold under \textsf{PFA} and we have proved that \textsf{PFA} implies no MAD family is almost-normal, we get the following:

\begin{cor}{(\textsf{PFA})}
Completely separable, strongly $(\aleph_0,<\mathfrak{c})$-separated MAD families which are not almost-normal exist generically.
\end{cor}

In particular, strongly $(\aleph_0,<\mathfrak{c})$-separated AD families which contain Luzin families (and hence are not potentially normal) exist generically.
We do not know if this is always the case, or at least, it follows from \textsf{MA}.

\begin{que}
Is it consistent that strongly $\aleph_0$-separated (or strongly $(\aleph_0,<\mathfrak{c})$-separated) AD families are potentially normal? Is it consistent with \textsf{MA}?
\end{que}

\subsection*{Acknowledgements}
The author wishes to thank Vinicius de Oliveira Rodrigues, Victor dos Santos Ronchim and Sergio Garc\'ia Balan for several stimulating discussions, and also Michael Hru\v{s}\'ak for some suggestions, in particular, for pointing out Theorem \ref{pfa}.

\bibliographystyle{plain}
\bibliography{Madness.bib}

\end{document}